\documentclass[11pt]{article}
\usepackage{fullpage,amsmath,amsthm}

\newtheorem{theorem}{Theorem}[section]

\newtheorem{proposition}[theorem]{Proposition}

\theoremstyle{definition}

\newcommand{\e}{\epsilon}

\newcommand{\angles}[1]{\langle #1 \rangle}

\newcommand{\poly}{\ensuremath{\mathrm{poly}}}

\title{A Simple Explicit Construction of an $n^{\Tilde{O}(\log n)}$-Ramsey Graph}

\author{Boaz Barak\thanks{Department of Computer Science, Princeton University
\texttt{boaz@cs.princeton.edu}.}}

\begin{document}

\maketitle

\begin{abstract}
We show a simple explicit construction of an
$2^{\Tilde{O}(\sqrt{\log n})}$ Ramsey graph. That is, we provide a
$\poly(n)$-time algorithm to output the adjacency matrix of an
undirected $n$-vertex graph with no clique or independent set of
size $2^{\e \sqrt{\log n}\log\log n}$ for every $\e>0$.

Our construction has the very serious disadvantage over the
well-known construction of Frankl and Wilson \cite{FranklWi81} that
it is only explicit and not very explicit, in the sense that we do
\emph{not} provide a poly-logarithmic time algorithm to compute the
neighborhood relation. The main advantage of this construction is
its extreme simplicity. It is also somewhat surprising that even
though we use a completely different approach we get a bound which
essentially equals the bound of \cite{FranklWi81}. This construction
is quite simple and was obtained independently by others as
well\footnote{P.~Pudlak, personal communications, July 2004.} but as
far as we know has not been published elsewhere.
\end{abstract}

\section{The Construction}

As mentioned above, we prove the following proposition:

\begin{proposition} Let $\e>0$ be some constant. There is a
polynomial-time  algorithm $A$ that on input $1^n$ outputs the
adjacency matrix for a graph $H$ on $n$ vertices with no clique or
independent set of size $2^{\e \sqrt{\log n}\log\log n}$.
\end{proposition}
\begin{proof} We will need to recall the notion of the Abbott
product of two graphs: if $G=(V_G,E_G)$ and $H=(V_H,E_H)$ are
graphs, then the \emph{Abbott product} of $G$ and $H$, denoted by $G
\otimes H$ is the graph with vertex set $V_H \times V_H$ (where
$\times$ denotes Cartesian product) and where
$\angles{(u,v),(u',v')}$ is an edge if either $(u,u')$ is an edge in
$G$ or $u=u'$ and $(v,v')$ is an edge in $H$. One can think of $G
\otimes H$ as obtained by replacing each node of $G$ with an entire
copy of $H$ (both vertices and edges), where each two different
copies of $H$ have either all the edges between them or none of the
edges between them, depending on whether the corresponding vertices
in $G$ are neighbors. We let $G^l$ denote $G \otimes G \otimes
\cdots \otimes G$ ($l$ times).

We let $\omega(G)$ be the clique number of $G$ (i.e., the size of
the largest clique in $G$) and $\alpha(G)$ be the independence
number of $G$ (i.e., the size of the largest independent set in
$G$). The basic fact we need about the Abbot product is that
$\omega(G \otimes H) = \omega(G) \cdot \omega(H)$ and $\alpha(G
\otimes H) = \alpha(G) \cdot \alpha(H)$.

We can now specify our construction. Given $\e>0$, the algorithm
$A'$ will choose a constant $c>1$ (the exact choice of $c$ will be
specified later), and let $k=2^{c\sqrt{\log n}}$ and using
$k^{O(\log k)}=n^{O(1)}$ running time construct a graph $G$ on $k$
vertices such that $\omega(G),\alpha(G) < 3\log k$.

Constructing such a graph can be done using well-known techniques:
as a first observation note that in time $k^{O(\log k)}$ it can be
\emph{verified} that a graph $G$ satisfies
$\omega(G),\alpha(G)<3\log k$. Thus, it is enough to show an
explicit family of  $k^{O(log k^2)}$ graphs, where one of which
satisfies this condition. A graph can be  represented as a string of
length $\binom{k}{2}$. We claim that if we choose this string from a
sample space that is $2^{-5\log^2k}$-close to being $5\log^k$-wise
independent then with high probability the graph will satisfy
$\omega(G),\alpha(G)<2\log k$, once we prove this then we'll be done
since explicit sample spaces with cardinality $k^{O(\log k)}$ were
given by Naor and Naor \cite{NaorNa93}. However, this follows by the
same reason that a random graph satisfies this property: that every
set of $4.5\log^2 k$ edges has probability at most $2\cdot
2^{-4.5\log^2 k} + 2^{-5\log^2 k} \ll 1/\binom{k}{3\log k}$ to be
identically zero or identically one.\footnote{Another approach that
may work is derandomization using the method of conditional
expectations.}

Now, the algorithm will compute the graph $H =
G^{\tfrac{1}{c}\sqrt{\log n}}$. This graph has $n$ vertices, but
\[ \omega(H),\alpha(H) < (3\log k)^{\tfrac{1}{c}{\sqrt{\log n}}} =
(3c\sqrt{\log n})^{\tfrac{1}{c}\sqrt{\log n}} < 2^{\tfrac{\log c +
2}{2c}\log\log n \sqrt{\log n}} \] we choose $c$ large enough such
that the constant expression in the exponent will be smaller than
$\e$.
\end{proof}


\begin{thebibliography}{LRVW03}

\bibitem[FW81]{FranklWi81}
P.~Frankl and R.~M.~Wilson.
\newblock Intersection theorems with geometric consequences.
\newblock {\em Combinatorica.}, 4(1):357--368,  1981.


\bibitem[NN93]{NaorNa93}
 J.~Naor and  M.~Naor
\newblock Small-Bias Probability Spaces: Efficient Constructions and Applications
\newblock {\em SIAM J. Comput.}, 22(4): 838-856 ,1993.

\end{thebibliography}
\end{document}